\documentclass[p4paper, 12pt]{article}
\usepackage{epsfig, color, epstopdf}
\usepackage{amssymb,amsmath,amsthm,latexsym,lineno}
\usepackage{graphicx,graphics,amsfonts,indentfirst}
\usepackage{mathrsfs,eufrak, bbm, dsfont,yfonts}
\usepackage{url,tikz,yhmath,appendix,booktabs}
\usepackage{multirow}
\usepackage{kantlipsum}

\allowdisplaybreaks
\usepackage{setspace}

\newtheorem{theo}{Theorem}
\newtheorem{lem}[theo]{Lemma}
\newtheorem{pro}[theo]{Proposition}
\newtheorem{prob}{Problem}
\newtheorem{exa}[theo]{Example}
\newtheorem{con}{Conjecture}
\newtheorem{cor}[theo]{Corollary}

\newcommand \equ[2]
{
	\begin{equation}\label{#1}
		#2
	\end{equation}
}

\newcommand \eqn[2]
{
	\begin{eqnarray}\label{#1}
		#2
	\end{eqnarray}
}

\newcommand \them[2]
{\begin{theo}
		\label{#1} #2
	\end{theo}
}

\newcommand \lemm[2]
{\begin{lem}
		\label{#1} #2
	\end{lem}
}

\newcommand \corr[2]
{\begin{cor}
		\label{#1} #2
	\end{cor}
}

\newcommand \conj[2]
{\begin{con}
		\label{#1} #2
	\end{con}
}

\newcounter{countcase}

\newcounter{countclaim}

\def \proof {\noindent {{\it Proof}}.\setcounter{countcase}{0} \setcounter{clm}{0}}

\newcommand{\proofend}{{\hfill$\Box$}}

\def \N {{\mathbb N}}

\def \broken {{\mathscr {B}}}
\def \nbroken {{\mathscr {N}}\hspace{-0.15 cm}{\mathscr {B}}}
\def \nbrokenf {{\mathscr {N}}\hspace{-0.15 cm}{\mathscr {B}}\hspace{-0.05 cm}{\mathscr {F}}}

\newcommand{\rebibitem}[1] 
{
	\bibitem{#1} 
}

\begin{document}
	
	\baselineskip 0.6 cm
	
\title{
	An improved lower bound of $P(G, L)-P(G,k)$ for $k$-assignments $L$
	}

\author{
		Fengming Dong\thanks{
			Corresponding Author. 
			Email: fengming.dong@nie.edu.sg and donggraph@163.com.} ~and  Meiqiao Zhang\thanks{Email: nie21.zm@e.ntu.edu.sg and 
								meiqiaozhang95@163.com.}
		\\
		National Institute of Education,
		Nanyang Technological University, Singapore
	}
	\date{}
	
	\maketitle{}

	\begin{abstract}
Let $G=(V,E)$ be a simple graph
with $n$ vertices and $m$ edges,
$P(G,k)$ be 
the chromatic polynomial of $G$, 
and $P(G,L)$ be the number of 
$L$-colorings of $G$ for any 
$k$-assignment $L$.
In this article, 
we show that when $k\ge m-1\ge 3$, 
$P(G,L)-P(G,k)$ is bounded below by 
$
\left ( (k-m+1)k^{n-3}
+\frac{(k-m+3)c}{3} k^{n-5}
\right )
\sum\limits_{uv\in E}|L(u)\setminus L(v)|$,
where $c\ge \frac{(m-1)(m-3)}{8}$,
and in particular, if $G$ is $K_3$-free,
then
$c\ge {m-2\choose 2}+2\sqrt m-3$.
Consequently, $P(G,L)\ge P(G,k)$ whenever $k\ge m-1$. 
	\end{abstract} 
	



\section{Introduction}
	\label{sec1}

For any graph $G$, let $V(G)$ 
and $E(G)$ be the vertex set and edge 
set of $G$ respectively.
Let $\N$ be the set of positive integers,
and any $k\in \N$, 
let $[k]:=\{1,\dots, k\}$.
A {\it proper  $k$-coloring} of $G$ 
is a map $f:V(G)\rightarrow [k]$ 
such that $f(u)\ne f(v)$ for 
each pair of adjacent vertices $u$ 
and $v$ in $G$. 
Let $P(G,k)$ denote the number of proper $k$-colorings of $G$. 
Introduced by Birkhoff~\cite{birk}
in 1912,  
$P(G,k)$ is called the {\it chromatic polynomial} of $G$.
More details on $P(G,k)$ can be found
in \cite{birk, birk2, dong1, dong0, Jackson2015, rea1, rea2}.

The notion of list-coloring was introduced independently 
by Vizing~\cite{viz} 
and by Erd\H{o}s, Rubin and Taylor~\cite{erdos}.
A map $L:V(G)\rightarrow 2^{\N}$ 
is  called an {\it assignment} of $G$. 
For any $k\in \N$,
a $k$-{\it assignment} of $G$ 
is an assignment $L$ of $G$ with  
$|L(v)|=k$ for all $v\in V(G)$. 
Given any assignment $L$ of $G$,
an {\it $L$-coloring} of $G$ is 
a map $f:V(G)\rightarrow \N$ with the property 
that $f(v)\in L(v)$ for each $v\in V(G)$
and $f(u)\ne f(v)$ for 
each pair of adjacent vertices $u$
and $v$ in $G$.
Let $P(G,L)$ denote the number of $L$-colorings of $G$.
For any $k\in \N$, 
let $P_l(G,k)$ be the minimum value of 
$P(G,L)$ among all $k$-assignments $L$ of $G$.
Introduced by Kostochka and Sidorenko in 1990s, 
 $P_l(G,k)$ is called {\it the list-color function} of $G$.
More details on $P_l(G,k)$ can be 
found in \cite{Thom}.

It is known that 
$P(G,k)$ is a polynomial in $k$ of degree $|V(G)|$
(see Theorem~\ref{th2-1}).
However, due to Donner~\cite{Donner},
$P_l(G,k)$ is in general not a polynomial in $k$.
By the definitions of $P(G,k)$ and $P_l(G,k)$, 
$P_l(G,k)\le P(G,k)$ holds for every 
$k\in \N$.
Clearly, $P_l(G,k)=P(G,k)$
does not hold for some graphs $G$ 
and some numbers $k\in \N$. 
For example, 
$P(G,2)\ge 2$ holds for each bipartite graph $G$,
but $P_l(G,2)=0$ as long as $G$ contains $K_{2,4}$ as a subgraph.
On the other hand, it is not difficult to verify that 
$P(G,k)=P_l(G,k)$ holds for 
any chordal graph $G$ and $k\in \N$
(see \cite{kos}).
From the big picture, 
for any simple graph $G$, 
Donner~\cite{Donner} 
showed that $P(G,k)=P_l(G,k)$ holds
when $k$ is sufficiently large,
answering a problem proposed by  Kostochka and Sidorenko \cite{kos},
and Thomassen~\cite{Thom} proved
that $P(G,k)=P_l(G,k)$ 
when $k>|V(G)|^{10}$.
In 2017, Wang, Qian and Yan~\cite{wangwei} significantly 
improved this result by showing that 
$P(G,k)=P_l(G,k)$ holds for any $k\in \N$ with 
$k\ge \frac{(m-1)}{\ln (1+\sqrt 2)}
\approx 1.135(m-1)$,
where $m$ is the number of edges in $G$.

In this article, we will establish 
a lower bound of $P(G,L)-P(G,k)$ for 
an arbitrary $k$-assignment $L$
with $k\ge m-1$.
Obviously, $P(G,L)=P(G,k)$ holds whenever $L(u)=L(v)$ for every edge $uv$ in $G$.
This article shows how large the gap between $P(G,L)$ 
and $P(G,k)$ can be when 
$L(u)\ne L(v)$ for some edge $uv$ in $G$.


\them{th1-1}
{
Let $G=(V,E)$ be a simple graph with $n$ vertices and $m\ (\ge 4)$ edges.
Then, for any $k$-assignment $L$ of $G$
with $k\ge m-1$,
\equ{main1}
{
P(G,L)-P(G,k)\ge  
\left ( 
(k-m+1)k^{n-3}
+\frac{(k-m+3)c}{3}k^{n-5}
\right )
 \sum_{uv\in E}|L(u)\setminus L(v)|,
}
where $c\ge \frac{(m-1)(m-3)}{8}$,
and particularly, 
when $G$ is $K_3$-free, 
 $c
\ge {m-2\choose 2}+2\sqrt m-3$.
}

Note that any graph with less than $4$ edges is a chordal graph.
Thus, the following conclusion follows from Theorem~\ref{th1-1} directly.

\corr{cor1-1}
{
For any simple graph $G$ with $m$ edges, $P_l(G,k)=P(G,k)$ holds for each 
$k\in \N$ with $k\ge m-1$. 
}

Let $G=(V,E)$ be a simple graph with $n$ vertices and $m$ edges and  
$\eta$ be a fixed bijection 
from $E$ to $[m]$.
A {\it broken cycle} of $G$ (with respect to $\eta$) 
is a path $B=v_1v_2\dots v_r$ of $G$,
where $r\ge 3$,
such that $v_1v_r\in E$ 
and $\eta(v_1v_r)<\eta(v_iv_{i+1})$ for each $i=1,2,\dots,r-1$.
Let $\broken(G)$ be the collection
of edge sets $E(B)$ over all broken cycles $B$ of $G$, and 
let $\nbroken(G)$ be 
the set  of subsets $A$ of $E$ 
that is broken-cycle free with 
respect to $\eta$
(i.e., $E_0\not\subseteq A$ 
for each $E_0\in  \broken(G)$).
Obviously, for each $A\in \nbroken(G)$,
the spanning subgraph $(V,A)$ has no cycles, implying that 
$0\le |A|\le n-1$.
For each $i$ with $0\le i\le n-1$,
let $\nbroken_i(G)$ be the set of $A\in \nbroken(G)$ with $|A|=i$.
 
For any $e\in E$ and $1\le i\le n-1$, let $\nbroken_i(G,e)$ 
be the set of $A\in \nbroken_i(G)$ 
with $e\in A$.
Note that $|\nbroken_i(G,e)|$ depends on $\eta$ although $\eta$ is not included in the notation. 
For example, if $G$ is $K_3$, then 
$|\nbroken_2(G,e)|$ is either $1$ or $2$. 
Let $Q_{\eta}(G,e,x)$ denote the polynomial defined below:
\equ{eq2-0}
{
	Q_{\eta}(G,e,x):
	=\sum_{1\le i\le n-1
		\atop i\ odd}\frac{|\nbroken_i(G,e)|}{i} x^{n-i}
	-\sum_{2\le i\le n-1
		\atop i\ even}|\nbroken_i(G,e)| x^{n-i}.
}
For any  $e\in E$, 
let $G/e$ denote the simple graph 
obtained from $G$ by contracting $e$ 
and deleting 
all but one of the multiple edges, if they arise.
Thus, $|E(G/e)|=m-1-t$, 
where $t$ is the number of $3$-cycles 
in $G$ containing $e$.

In Section~\ref{sec2}, 
we show that if $x\ge m-1$ and $n\ge 4$,  then 
$$
Q_{\eta}(G,e,x)\ge
(x-m+1)x^{n-2}+
\frac {(x-m+3)|\nbroken_2(G/e)|}{3}x^{n-4}.
$$
Then, in Section~\ref{sec30}, 
we find a lower bound of $|\nbroken_2(G/e)|$ in terms of $m$.
In Section~\ref{sec3},
we prove that 
$P(G,L)-P(G,k)$
is bounded below by
$\frac 1k \sum\limits_{uv\in E}(|L(u)\setminus L(v)|Q_{\eta}(G,uv,k))$
for any $k$-assignment $L$ 
of $G$ with $k\ge m-1$.
Theorem~\ref{th1-1}
then follows immediately.
Finally, in Section~\ref{last sec}, 
we propose two conjectures 
studying the relation between 
$P_l(G,k)$ and $P(G,k)$. 

\section{A lower bound of 
	$Q_{\eta}(G,e,x)$}
\label{sec2}

In this section, we always assume that 
$G=(V,E)$ is a simple graph with $n$ vertices and $m$ edges and  
$\eta$ is a fixed bijection 
from $E$ to $[m]$.
Due to Whitney~\cite{Whi1932},
the coefficients of $P(G,x)$ can be expressed in terms of the sizes of $\nbroken_i(G)$'s.

\begin{theo}[\cite{Whi1932}]
\label{th2-1}
$P(G,x)$ can be expressed as 
$P(G,x)=\sum\limits_{i=0}^{n-1} (-1)^i 
|\nbroken_i(G)| x^{n-i}$.
\end{theo}

In this section, we shall find a lower bound of $Q_{\eta}(G,e,x)$ 
for any edge $e$ under the condition 
$x\ge m-1$. 
By the definition of $\nbroken_i(G,e)$, 
we first have the following relation between $|\nbroken_i(G,e)|$ and $|\nbroken_{i+1}(G,e)|$.

\lemm{le2-1}
{
For any $e\in E$
and $ i\in [n-2]$,
$i|\nbroken_{i+1}(G,e)|\le (m-i)|\nbroken_i(G,e)|$.
}

\proof 
When $i\ge m$, the inequality is trivial, as both sides are $0$.
Now assume that $1\le i\le m-1$. 
Lemma~\ref{le2-1} then follows directly from the following facts:
\begin{enumerate}
	\item for each $A\in \nbroken_{i+1}(G,e)$ and $e'\in A\setminus \{e\}$, $A\setminus \{e'\}\in \nbroken_i(G,e)$; and 
	\item for each $A'\in \nbroken_{i}(G,e)$, there are at most 
	$m-i$ edges $e'$ in $E\setminus A'$ such that 
	$A'\cup \{e'\}\in \nbroken_{i+1}(G, e)$.
	\proofend
\end{enumerate}

We can now apply Lemma~\ref{le2-1}
to find a lower bound of $Q_{\eta}(G,e,x)$.

\them{the2-2}
{
Assume that $n\ge 3$.
For any edge $e$ in $G$ and $x\ge 0$, 
\equ{eq2-3}
{
Q_{\eta}(G,e,x)
\ge \sum_{1\le i\le n-1
	\atop i\ odd}\frac{|\nbroken_i(G,e)|}{i} 
(x-m+i) x^{n-i-1}.
}
In particular, if $n$ is even, then,
\equ{eq2-31}
{
Q_{\eta}(G,e,x)
\ge \sum_{1\le i\le n-3
		\atop i\ odd}
	\frac{|\nbroken_i(G,e)|}{i} 
	(x-m+i) x^{n-i-1}+
\frac{|\nbroken_{n-1}(G,e)|}{n-1}x.
}
}

\proof 
By Lemma~\ref{le2-1},
for any $i\in [n-2]$, 
as $x\ge 0$, 
\eqn{eq2-4}
{
\frac {|\nbroken_i(G,e)|}{i}x^{n-i}-
|\nbroken_{i+1}(G,e)|x^{n-i-1}
&\ge &\frac{|\nbroken_i(G,e)|}{i}x^{n-i}-\frac{(m-i)|\nbroken_i(G,e)|}{i}x^{n-i-1}
\nonumber \\
&=&\frac{|\nbroken_i(G,e)|}{i}(x-m+i)x^{n-i-1}.
}
By the definition of $Q_{\eta}(G,e,x)$,
the result follows from 
(\ref{eq2-4}).
\proofend

For any edge $e$ in $G$, let 
$\eta|_{E(G/e)}$ be the restriction of $\eta$  to the  edge set of $G/e$, and let $\nbroken_{j}(G/e)$ 
be the set of $A\subseteq E(G/e)$ 
with $|A|=j$ 
such that $A$ is broken-cycle free 
with respect to $\eta|_{E(G/e)}$.
In the following, we will show that 
$|\nbroken_i(G,e)|$ is bounded below
by $|\nbroken_{i-1}(G/e)|$.

\lemm{le2-2}
{
For any $e\in E$ and $i\in [n-1]$, 
$
|\nbroken_i(G,e)|\ge |\nbroken_{i-1}(G/e)|.
$
}
\proof
It suffices to show that $A\cup \{e\}\in \nbroken_i(G,e)$ for each $A\in \nbroken_{i-1}(G/e)$.

Suppose that $A\cup \{e\}\not\in \nbroken_i(G,e)$. 
Then, there exists $B\in \broken(G)$ 
with $B\subseteq A\cup\{e\}$.
As $A\in \nbroken_{i-1}(G/e)$, $B\not\subseteq A$, which implies that $e\in B$ and $B\setminus \{e\}\subseteq A$. 
However, $B\in \broken(G)$ implies that 
$B\setminus \{e\}\in \broken(G/e)$,
a contradiction to the assumption that $A\in \nbroken_{i-1}(G/e)$.

Hence Lemma~\ref{le2-2} follows.
\proofend

Combining Theorem~\ref{the2-2} and Lemma~\ref{le2-2}, we obtain a lower bound of $Q_{\eta}(G,e,x)$
in terms of $|\nbroken_2(G/e)|$ and $x$.

\corr{cor2-1}
{
For any $e\in E$ and real number $x$ with $x\ge m-1$, 
if $n\ge 4$,  then 
\equ{eq2-6}
{
	Q_{\eta}(G,e,x)
	\ge 
	(x-m+1)x^{n-2}+
	\frac{(x-m+3)|\nbroken_2(G/e)|}{3} x^{n-4}.
}
}

\section{Lower bounds of $|\nbroken_2(G/e)|$}
\label{sec30}

In this section, we 
still assume that $G=(V,E)$ is
a simple graph with $|V|=n$ and 
$|E|=m$,
and we shall find a lower bound of $|\nbroken_2(G/e)|$ in terms of
$m$ for an arbitrary edge $e$ in $G$.

Given any simple graph $H$, 
by the definition of $|\nbroken_2(H)|$ 
or Corollary 2.3.1 in \cite{dong0}, 
$|\nbroken_2(H)|$ has the following expression:
\equ{eq2-30}
{
	|\nbroken_2(H)|=
	\binom{|E(H)|}{2}
	-\triangle(H),
}
where $\triangle(H)$ is the number of 
$3$-cycles in $H$.

First consider the special case that $G$ is $K_3$-free.
Let $c_4(G)$ be the minimum integer $r$ such that each edge 
$e$ in $G$ is 
contained in at most $r$ $4$-cycles of $G$.
For any $u\in V$, let $N_G(u)$ denote the set of vertices in $G$ adjacent to $u$, and let $d_G(u)=|N_G(u)|$.

\lemm{le2-30}
{For any $e\in E$,
if $G$ is $K_3$-free and $m\ge 3$, then 
\equ{eq2-30-1}
{
|\nbroken_2(G/e)|\ge {m-1\choose 2}-c_4(G)
\ge {m-2\choose 2}+2\sqrt m-3.
}
}

\proof 
As $G$ is $K_3$-free, 
then $G/e$ has exactly $m-1$ edges and  at most $c_4(G)$ 
$3$-cycles.
 Thus, applying (\ref{eq2-30})  implies  that 
 $|\nbroken_2(G/e)|\ge {m-1\choose 2}-c_4(G)$ for any edge $e\in E$.

Note that 
 ${m-1\choose 2}-{m-2\choose 2}-2\sqrt m+3=(\sqrt m-1)^2$.
 Thus, it remains to show that 
$c_4(G)\le (\sqrt m-1)^2$.
It suffices to show that
for each edge $e'$ in $G$, 
the number of $4$-cycles in $G$ 
containing $e'$, denoted by $c_4(e')$, 
is at most $(\sqrt m-1)^2$.
Let $e'=uv\in E$,
$N'(u):=N_G(u)\setminus \{v\}=\{u_1,u_2,\dots,u_p\}$
and $N'(v):=N_G(v)\setminus \{u\}=\{v_1,v_2,\dots,v_q\}$.
As $G$ is $K_3$-free, 
$N'(u)\cap N'(v)=\emptyset$.
If $p=0$ or $q=0$, then 
$c_4(e')=0< (\sqrt m-1)^2$.
Now, assume $p\ge 1$ and $q\ge 1$.
Clearly, $c_4(e')$ 
is equal to the size of the edge set 
$E_G(N'(u),N'(v)):=\{u_iv_j\in E: i\in [p], j\in [q]\}$. 
Thus, 
\equ{eq2-32}
{
	c_4(e')=|E_G(N'(u),N'(v))|
	\le m-1-p-q,
}
implying that $p+q\le m-1-c_4(e')$,
and therefore 
\equ{eq2-33}
{
	c_4(e')=|E_G(N'(u),N'(v))|\le pq\le \frac 14 (p+q)^2
	\le 
	\frac 14 (m-1-c_4(e'))^2.
}
Solving the inequality 
$c_4(e')\le \frac 14(m-1-c_4(e'))^2$ 
with the condition $c_4(e')
<m$ gives that $c_4(e')\le (\sqrt m-1)^2$.
Hence $c_4(G)=\max\limits_{e'\in E} c_4(e')\le (\sqrt m-1)^2$. The result holds. 
\proofend

Now we are going to find a lower bound
of ${|\nbroken_2(G/e)|}$
in terms of $m$ for 
any edge $e$ in $G$.
We shall apply the following theorem obtained by Fisher in~\cite{Fish1989}.

\begin{theo}[\cite{Fish1989}]
\label{th2-4}
For any simple graph $H$, 
$\triangle(H)\le \frac 16 |E(H)|(\sqrt{8|E(H)|+1}-3)$.
\end{theo}

	By applying Theorem~\ref{th2-4},
we can find an upper bound of 
$\triangle(H)$ in terms of $|E(H)|$
and $t$, where $t$ is any number 
not larger than the maximum degree of $H$.

\lemm{le2-3-1}
{
For any simple graph $H$, 
if the maximum degree of $H$ is at least $t$, then 
\eqn{eq3-60-0}
{
\triangle(H)\le 
\frac{|E(H)|-t}6
\left (3
+ \sqrt{8(|E(H)|-t)+1}\right).
}
}

\proof 
Let $w$ be a vertex in $H$ 
with $d_H(w)=s\ge t$.
Let $H_0$ be the subgraph of $H$ 
induced by $N_H(w)$, and 
let $H-w$ be the subgraph of $H$ induced by $V(H)\setminus \{w\}$.
Then $|E(H_0)|\le |E(H)|-s$ and $|E(H-w)|=|E(H)|-s$.
Then,  
\eqn{eq3-60}
{
\triangle(H)
&=&|E(H_0)|
+\triangle(H-w)
\nonumber \\
&\le & |E(H)|-s
+ \frac 16 (|E(H)|-s)
\left (\sqrt{8(|E(H)|-s)+1}-3\right )
\nonumber \\
&=&
\frac{|E(H)|-s}6
\left (3
+ \sqrt{8(|E(H)|-s)+1}\right) ,
}
where the penultimate expression follows from Theorem~\ref{th2-4}.
As $s\ge t$, 
the lemma holds.
\proofend

\lemm{le2-3}
{If $m\ge 4$, then 
	for any $e\in E$,
	$
	|\nbroken_2(G/e)|
	\ge \frac{(m-1)(m-3)}{8}.
	$
}

\proof 
Let $e$ be any edge in $G$ and 
let $t$ be the number of $3$-cycles in $G$ containing $e$. 
Then $m\ge 2t+1$ and $|E(G/e)|=m-t-1$.
By (\ref{eq2-30})
and Theorem~\ref{th2-4}, 
$|\nbroken_2(G/e)|\ge g(t,m)$, where 
\eqn{eq3-61}
{
	g(t,m)
	:&=& {m-t-1\choose 2}-
	\frac{(m-t-1)}6
	\left (\sqrt{8(m-t-1)+1}-3\right )
	\nonumber \\
	&=&\frac{(m-t-1)^2}{2}-
	\frac{(m-t-1)}{6}\sqrt{8(m-t-1)+1}.
}
Note that $f(x):=\frac 12 x^2-\frac x6 \sqrt{8x+1}$ is strictly increasing for $x\ge 1$, implying that $f(m-1)\ge f(m-2)$. 
	Since $g(t,m)=f(m-1-t)$, 
it is routine to verify that
when $m\ge 4$, 
\eqn{eq3-62}
{
g(0,m)>g(1,m)
=\frac{(m - 2)^2}2 - \frac{m-2}6
\sqrt{8m - 15}
>\frac{(m-1)(m-3)}8.
}
It remains to consider the case $t\ge 2$.
Note that $|E(G/e)|=m-t-1$
and the vertex in $G/e$ produced 
after contracting $e$ is of degree at least $t$.
As $m\ge 2t+1$,
by (\ref{eq2-30}) and  Lemma~\ref{le2-3-1}, 
\eqn{eq3-64}
{
|\nbroken_2(G/e)|
&\ge & {m-t-1\choose 2}
-\frac{(m-2t-1)}6
\left (3+\sqrt{8(m-2t-1)+1}\right )
\nonumber \\
&=&\frac{(m-1)(m-3)}{8}
+	\frac{m-2t-1}{24}
\left (
9m-6t-27-4\sqrt{8(m-2t-1)+1}
\right )
\nonumber \\
&\ge &\frac{(m-1)(m-3)}{8}.
}	
Hence Lemma~\ref{le2-3} holds.
\proofend

By Corollary~\ref{cor2-1} and
Lemmas~\ref{le2-30} and~\ref{le2-3},
the following conclusion holds.

\them{th2-2}
{
	For any $e\in E$ and real number $x$ with $x\ge m-1\ge 3$,  
$Q_{\eta}(G,e,x)$ is bounded below by 
$(x-m+1)x^{n-2}+
\frac{(x-m+3)c}{3} x^{n-4}$,
where $c\ge \frac{(m-1)(m-3)}{8}$,
and in particular, if $G$ is $K_3$-free,
then $c
\ge {m-2\choose 2}+2\sqrt m-3$.
}

\section{Proving Theorem~\ref{th1-1}}\label{sec3}

In this section, we always assume that
$G=(V,E)$ is a simple graph with $n$ vertices and $m$ edges, 
$\eta$ is a fixed bijection 
from $E$ to $[m]$, 
and $L$ is a $k$-assignment of $G$,
where $k\ge 2$.

For any integer $i$ 
with $0\le i\le n-1$,
let $\nbrokenf_i(G)$ be the set 
of spanning forests $F=(V,A)$ of $G$ 
with $A\in \nbroken_i(G)$.
Clearly, each $F\in \nbrokenf_i(G)$ has
exactly $n-i$ components.
We can represent $F$ by the set  $\{T_1,T_2,\dots,T_{n-i}\}$,
where $T_1,T_2,\dots,T_{n-i}$
are the components of $F$.

For any subgraph $H$ of $G$, define 
$\beta(H)=\big | 
\bigcap\limits_{v\in V(H)}L(v)
\big |$.
By applying the inclusion-exclusion principle, it can be proved that 
\equ{eq3-11}
{
P(G,L)=\sum_{i=0}^{n-1} 
(-1)^i 
\sum_{\{T_1,\dots,T_{n-i}\}\in \nbrokenf_i(G)} \prod_{j=1}^{n-i}
\beta(T_j).
}
By Theorem~\ref{th2-1} and (\ref{eq3-11}), we have 
\equ{eq3-12}
{
	P(G,L)-P(G,k)=
	\sum_{i=1}^{n-1} 
	(-1)^i 
	\sum_{\{T_1,\dots,T_{n-i}\}\in \nbrokenf_i(G)} 
	\left (\prod_{j=1}^{n-i}
	\beta(T_j)-k^{n-i}\right ).
}
For any edge $e=uv$ in $G$, let 
$\alpha(e)=|L(u)\setminus L(v)|$.
For any $F=\{T_1,\dots,T_{n-i}\}\in \nbrokenf_i(G)$,
a lower bound for $\prod\limits_{j=1}^{n-i}
\beta(T_j)-k^{n-i}$
 was obtained in \cite{wangwei},
 as stated below.

\begin{lem}[\cite{wangwei}]
	\label{le3-3}
For any $i\in [n-1]$ and $F=\{T_1,\dots,T_{n-i}\}\in \nbrokenf_i(G)$,
\equ{eq3-13}
{
\prod_{j=1}^{n-i}
\beta(T_j)-k^{n-i}
\ge -k^{n-i-1} \sum_{e\in E(F)}\alpha(e).
}
\end{lem}

We are now going to establish an upper bound for $\prod\limits_{j=1}^{n-i}
\beta(T_j)-k^{n-i}$.
We first introduce the following result.

\lemm{le3-1}
{
Let $d_1,d_2,\dots,d_r$ be 
any non-negative real numbers,
and $q_1,q_2,\dots,q_r$ be any positive 
real numbers,
where $r\ge 1$.
If $x\ge \max\limits_{1\le i\le r}d_i$, then 
\equ{eq3-1}
{
(x-d_1)(x-d_2)\cdots (x-d_r)
\le x^{r}-\frac {x^{r-1}}
{q_1+\cdots+q_r}\sum_{i=1}^r q_id_i.
}
}

\proof Assume that $d_1\ge d_2\ge \dots\ge d_r$. 
It is trivial to verify that 
$d_1\ge \frac 1{q_1+\cdots+q_r}\sum\limits_{i=1}^r q_id_i$.
As $0\le (x-d_i)\le x$ for all 
$2\le i\le r$, the result follows immediately. 
\proofend

\lemm{le3-4}
{For $i\in [n-1]$ and 
$F=\{T_1,\dots,T_{n-i}\}\in \nbrokenf_i(G)$,
	\equ{eq3-13'}
	{
		\prod_{j=1}^{n-i}
		\beta(T_j)-k^{n-i}
		\le -\frac{k^{n-i-1}}{i} \sum_{e\in E(F)}\alpha(e).
	}
}

\proof 
For each $T_j$ with $E(T_j)\ne \emptyset$, we have 
\equ{eq3-15}
{
\beta(T_j)
\le \min_{uv\in E(T_j)}|L(u)\cap L(v)|
\le k-\max_{e\in E(T_j)}\alpha(e)
\le k-\frac 1{|E(T_j)|}
\sum_{e\in E(T_j)}\alpha(e).
} 
Assume that $E(T_j)\ne \emptyset$ 
for each $j$ with $1\le j\le s$ 
while $E(T_j)= \emptyset$ for each $j$ with $s+1\le j\le n-i$.
As $|E(T_1)|+\cdots+|E(T_s)|=i$
and $k\ge \alpha(e)$ for each $e\in E(G)$, 
by (\ref{eq3-15})
and Lemma~\ref{le3-1}, 
\eqn{eq3-16}
{
\prod_{j=1}^{s} \beta(T_j)
&\le & 
	\prod_{j=1}^{s}
	\left (k-\frac 1{|E(T_j)|}
	\sum_{e\in E(T_j)}\alpha(e)\right )
	\nonumber \\
&\le & 
k^{s}-\frac {k^{s-1}}{|E(T_1)|+\cdots+|E(T_s)|}
\sum_{j=1}^s\sum_{e\in E(T_j)}\alpha(e)
	\nonumber \\
&= & 
k^{s}-\frac {k^{s-1}}{i}
\sum_{e\in E(F)}\alpha(e).
} 
As $\beta(T_j)=k$ for each $j$ with $s+1\le j\le n-i$, (\ref{eq3-13'}) follows.
\proofend

Recall that $Q_{\eta}(G,e,k)$ is the function defined in (\ref{eq2-0})
and $\nbroken_i(G,e)$ is the set 
of $A\in \nbroken_i(G)$ with $e\in A$.
We are now going to
find a lower bound of $P(G,L)-P(G,k)$
in terms of $\alpha(e)$ and $Q_{\eta}(G,e,k)$ for all edges $e$ in $G$.

\lemm{le3-5}
{$
P(G,L)-P(G,k)
\ge \frac 1k 
\sum\limits_{e\in E}
\left ( \alpha(e) Q_{\eta}(G,e,k)\right ).
$
}

\proof
By (\ref{eq3-12})
and applying 
Lemma~\ref{le3-3} for even $i$'s
and Lemma~\ref{le3-4} for odd $i$'s,
\eqn{eq3-17}
{
& &P(G,L)-P(G,k)
\nonumber \\\
&=&\sum_{i=1}^{n-1}
(-1)^i 
	\sum_{\{T_1,\dots,T_{n-i}\}\in \nbrokenf_i(G)} 
\left (\prod_{j=1}^{n-i}
\beta(T_j)-k^{n-i}\right )
\nonumber \\
&\ge &
\sum_{1\le i\le n-1 \atop i\ odd}
\frac{k^{n-i-1}}{i}
\sum_{F\in \nbrokenf_i(G)}\sum_{e\in E(F)}\alpha(e)
-\sum_{2\le i\le n-1 \atop i\ even}
k^{n-i-1}\sum_{F\in \nbrokenf_i(G)}\sum_{e\in E(F)}\alpha(e)
\nonumber \\
&=&
\sum_{1\le i\le n-1 \atop i\ odd}
\frac{k^{n-i-1}}{i}
\sum_{e\in E}
\alpha(e) |\nbroken_{i}(G,e)|
-\sum_{2\le i\le n-1 \atop i\ even}
k^{n-i-1}
\sum_{e\in E}
\alpha(e) |\nbroken_{i}(G,e)|
\nonumber \\
&=&\frac 1k
\sum_{e\in E}
\alpha(e) 
\left (
\sum_{1\le i\le n-1 \atop i\ odd}
\frac{|\nbroken_{i}(G,e)|}{i} 
k^{n-i}
-
\sum_{2\le i\le n-1 \atop i\ even}
|\nbroken_{i}(G,e)| k^{n-i} 
\right ).
}
By the definition of $Q_{\eta}(G, e, k)$, the result follows.
\proofend 

\vspace{0.2 cm}
We are now going to prove 
Theorem~\ref{th1-1}.

\noindent {\it Proof of Theorem~\ref{th1-1}}: 
As $m\ge 4$ and $k\ge m-1$,
 by Theorem~\ref{th2-2}
and Lemma~\ref{le3-5}, 
\eqn{eq3-20}
{
P(G,L)-P(G,k)
&\ge&    
\frac 1k \sum_{e\in E}(\alpha(e)
Q_{\eta}(G,e,k))
\nonumber \\
&\ge &
\left ( 
(k-m+1)k^{n-3}+\frac {(k-m+3)c}{3}  k^{n-5}
\right )
\sum_{e\in E} 
\alpha(e),
}
where 
$c\ge \frac{(m-1)(m-3)}{8}$,
and if $G$ is $K_3$-free,
then $c
\ge {m-2\choose 2}+2\sqrt m-3$.
\proofend 

\section{Concluding remarks
\label{last sec} }

\def \hyh {{\cal H}}

Given any simple graph $G$,  
the {\it list-chromatic number} 
of $G$, denoted by $\chi_l(G)$, 
is the minimum 
integer $r$ with $P_l(G,r)>0$,
and 
the {\it list-color function threshold} of $G$, denoted by $\tau(G)$, 
is the smallest integer $r \ge \chi(G)$ such that $P_l(G,k)= P(G,k)$ whenever 
$k\ge r$.
Obviously, $\tau(G)\ge \chi_l(G)\ge \chi(G)$.
By Corollary~\ref{cor1-1}, 
$\tau(G)\le |E(G)|-1$ when $|E(G)|\ge 4$.
The authors of this article also found some results on the upper bounds of
$\tau(\hyh)$  for
$r$-uniform hypergraphs $\hyh$ with $m$ edges, where $r \ge 3$.
In \cite{dong2},  
they showed that 
 $\tau(\hyh)$ is bounded above by $\min\{m-1, 0.6(m-1)+0.5\gamma(\hyh)\}$,
where $\gamma(\hyh)=\max_{e\in E(\hyh)}|E_{r-1}(e)|$
and $E_{r-1}(e)$ is the set of edges 
$e'$ in $\hyh$ with $|e \cap e'| = r -1$,
and in \cite{dong3},
they further showed that 
if $\rho(\hyh):=\min_{e,e'\in E(\hyh)}|e\setminus e'|\ge 2$ 
and $m\ge \rho(\hyh)^3/2+1$, 
then 
 $\tau(\hyh)\le \frac{2.4(m-1)}{\rho(\hyh)\log (m-1)}$.

Thomassen~\cite{Thom} asked 
if there exists a universal constant 
$\alpha$ such that $\tau(G)- \chi_l(G)\le \alpha$ holds for every simple graph $G$.  
Recently, Kaul et al \cite{Kaul2022}
gave a negative answer to Thomassen's question by showing that 
$\tau(K_{2,s})-\chi_l(K_{2,s})
\ge C\sqrt s$ for
a constant $C$ and 
all $s\in \N$ with $s\ge 16$.

We end this article with two conjectures.

\conj{con1}
{There exists a constant $c>0$ such that
$\tau(G)\le c|V(G)|$ 
or $\tau(G)\le c\Delta(G)$
for every simple graph $G$, 
where $\Delta(G)$ is the maximum degree of $G$.
}

\conj{con2}
{For any simple graph $G$, if $L$ is a 
$k$-assignment  of $G$ with $k\ge \tau(G)$ and 
$L(u)\ne L(v)$ for some edge $uv$ in $G$, 
then $P(G,L)>P(G,k)$.
}

Clearly, Conjecture~\ref{con2} holds for all chordal graphs. If this conjecture fails, then there exists a non-chordal graph $G$ 
such that $P(G,L)=P(G,k)$ for some 
$k$-assignment $L$, where $\tau(G)\le k\le |E(G)|-2$ and $L(u)\ne L(v)$ for some edge $uv$ in $G$.

\section*{Acknowledgement}

We thank the referees for their helpful comments,
and the second author would like to express her gratitude to National Institute of Education and Nanyang Technological University of Singapore for offering her Nanyang Technological University Research Scholarship during her Ph. D. study.

This research is supported by the Ministry of Education,
Singapore, under its Academic Research Tier 1 (RG19/22). Any opinions,
findings and conclusions or recommendations expressed in this
material are those of the authors and do not reflect the views of the
Ministry of Education, Singapore.

\end{document}